\documentstyle[12pt,A4]{article}

\newcommand{\bq}{\begin{equation}}
\newcommand{\eq}{\end{equation}}
\newcommand{\ba}{\begin{array}}
\newcommand{\ea}{\end{array}}
\newcommand{\ds}{\displaystyle}

\def\f2{\frac{1}{2}}
\def\nn{\mbox
{\sf
I\kern-.17em\rule[1.370ex]{.17em}{1pt}\kern-.17em\rule{.17em}{1pt}\kern-.12em{N}}}
\def\mm{\mbox
{\sf
I\kern-.17em\rule[1.370ex]{.17em}{1pt}\kern-.17em\rule{.17em}{1pt}\kern-.12em{M}}}

\title{\large\bf On the Distribution Function
of the Complexity of Finite Sequences}

\author{ J. Szczepa\'nski \\[-0.4cm]
{\small\it Polish Academy of Sciences,} \\[-0.4cm]
{\small\it  Institute of Fundamental
Technological Research}\\[-0.4cm]
{\small\it \'Swi\c{e}tokrzyska 21, 00-049 Warsaw, Poland}\\[-0.4cm]
{\small\it e-mail: jszczepa@ippt.gov.pl}\\[0.4cm]
}
\date{}
\begin{document}
\maketitle

{\small\bf {\it \bf Abstract} - Investigations
of complexity of sequences lead to important
applications such as effective data compression, testing of randomness,
discriminating between information sources and many others. In this
paper we establish formulas describing the distribution functions
of random variables representing the complexity of finite sequences
introduced by Lempel and Ziv in 1976. We show that the distribution
functions depend in an affine way on the probabilities of the so called
"exact" sequences.

\vspace{.2cm}

{\it\bf Keywords} : Complexity of sequence,
distribution
function, combinatorial problems, Lempel-Ziv parsing algorithms, randomness}

\vspace{.5cm}
\centerline{\sc I. Introduction}
The notion of complexity of a given sequence was first introduced
in papers by Kolmogorov [3] and Chaitin [1]. Kolmogorov proposed to
use the length of the shortest binary program which, when fed into
a given algorithm, will cause it to produce a specified sequence, as a
measure for the complexity of that sequence with respect to the given
algorithm. If the length of the program is large we can say that the
complexity of the sequence is large.

In 1976 Lempel and Ziv [4] proposed and explored another approach to
the problem of the complexity of a specific sequence. They linked the
complexity of a specific sequence to the gradual buildup of new
patterns along the given sequence. The complexity measure suggested by
them is related to the number of distinct phrases and the rate of
their occurence along the sequence. It reflects the behaviour of a
simple parsing algorithm whose task is to recognize newly encountered
phrases during its scanning of a given sequence. In a series of papers,
modifications of the Lempel-Ziv  parsing algorithm were proposed
in response to the needs of various applications. In general, in these
algorithms a new phrase is established as the shortest substring which
has not occurred previously, where the search for previous occurrences
may be restricted or generalized in the modified algorithms in various
ways, e.g.: by considering only a fixed number of preceding symbols
[8], by considering only complete previously established phrases
(Lempel-Ziv Incremental Parsing Algorithm [9]), by allowing a number
(not more than a fixed threshold) of previous occurrences of the phrase
(Generalized Lempel-Ziv Algorithm [6]) etc..

It turned out that investigations of sequence complexity
play an important role in universal data compression
schemes and their numerous applications such as efficient transmission
of data [8],[9], tests of randomness [16], discriminating between
information sources [2], [10], estimating the statistical model
of individual sequences [10] and many others.

In this paper we introduce the concept of exact sequences
i. e. sequences in which the last phrase of the sequence does not
occur in the past (precise formulation: Def. 3). We derive formulas
describing the distribution function of random variables representing
the complexity of finite sequences as defined by Lempel and Ziv in
1976. These formulas turn out to be of affine form with respect to
the probabilities of exact sequences.

\vspace{.5cm}
\centerline{\sc II. Lempel-Ziv complexity}

In this section we introduce the notation and recall basic
definitions [4].

Let \ ${\cal A}$ \ be a finite alphabet and let \ $\alpha =
|{\cal A}|$ \ denote  the size of the alphabet. Let \ ${\cal A}^{n} $
\ be the set of all sequences of length \ $n$ \ over \ ${\cal A}$ \
and let \ $S=s_1s_2\ldots s_n$ \ be an arbitrary element of \
${\cal A}^n\ .$ By \ $S(i,j)$ \ we denote the substrings \ $s_i
s_{i+1}\ldots s_j$ \ of \ $S$ \ when \ $ i \leq j$ \ and \ $S(i,j)
= \Lambda$ \ when \ $j < i\ .$ The partition
\bq\label{1}
H(S) = S(1, h_1) S(h_1 + 1,h_2)\ldots S(h_{m-1} + 1,n)
\eq
of \ $S$ \ such that for every \ $i ,\ S(h_{i-1} + 1,h_i -1)$ \ is a
substring of \ $S(1,h_i - 2)$ \ is called the history of \ $S$ \
and the \ $m$ \ strings \ $H_i(S) = S(h_{i-1} + 1, h_i)\ , \
i=1,2,\ldots,m$ \ where \ $h_0 = 0$ \ and \ $h_m = n\ ,$ are called
the components of the history. (Note that \ $h_1 = 1).$
Let \ $ c_H
(S)$ \ denote the number of components in a history \ $H(S) $ \ of
\ $S\ . $

{\it Definition 1}: The complexity \ $ c(S)$ \ of the
sequence \ $S$ \ is the number
\bq\label{2}
c(S) = \min \{ c_H (S) \}
\eq
where the minimum is over all histories of \ $S\ .$

{\it Definition 2}: The component \ $H_i(S)=S(h_{i-1}
+1,h_i)$ \ is called {\it exhaustive} if this string does not appear in
the string \ $S(1, h_i - 1)\ .$ A history of \ $S$ \ is called
{\it exhaustive} if each of its components, except possibly the last one,
is exhaustive.

It is easy to see that every sequence has a unique
exhaustive history, denoted by \ $H_E(S)\ .$
For instance, the exhaustive history
of the sequence
S=0011011101110110 is given by the following parsing of S :
0, 01, 10, 111, 0110110 where successive components are separated
by commas.

{\it Remark 1}: It was proved in [4] that \ $c(S) =c_E(S)\
, $ where \ $ c_E(S)$ \ is the number of components in \ $H_E (S)\
.$ Thus, below we shall use \ $c_E$ \ as the definition of complexity.

{\it Definition 3}: The sequence \ $S = s_1s_2\ldots s_n$
\ is called {\it exact} if the last string \ $S(h_{m-1} + 1,n)$ \ in its
exhaustive history \ $H_E(S) = S(1,h_1) S(h_1 + 1,h_2)\ldots
S(h_{m-1} + 1, n)$ \ does not occur as a substring \ $ S(i,j)$ \
(where \ $1\leq i \leq j \leq n - 1)$ in the sequence \ $S(1,n-1) =
s_1\ldots s_{n-1}\ .$

From now on we shall assume that for a fixed \ $n$ \ any element of \
$ {\cal A}^n$ \ is equi-probable, i.e. we assign the same
probability \ $\alpha^{-n}$ \ to each element of \ $ {\cal A}^n$ \
and
\bq\label{3}
P_n : 2^{{\cal A}^n}  \rightarrow [0,1]
\eq
denotes the probability in this sense. By $P_n(k)$ we denote
the probability of the event consisting of all sequences of length
 $n$ and complexity $k$  while $P_n(k_e)$ is the
probability of the event consisting of all exact sequences of length
$n$ and complexity~~$k.$


Under the above assumptions for every \ $n \in \nn $ we define the
random variable \ $C_n: {\cal A}^n \rightarrow \nn$ \ representing
the complexity:
\bq\label{4}
C_n(S): = c_E (S)
\eq
for every sequence \ $S \in {\cal A}^n\ .$

\vspace{.5cm}
\centerline{\sc III. The distribution function of \ $C_n$}

In this section we describe the distribution function of \ $C_n\ ,
\ n \in \nn\ .$ We prove the following

{\bf Theorem}: Under the above notation,
\bq\label{5}
P_{n+1} (C_{n+1} \leq k) = 1 - {\ds\sum\limits^n_{r=1}} P_r (k_e)
\eq
for every \ $n, k \in \nn\ .$

 {\it Proof:}
 We first express \ $P_{n+1} (k+1)$ \ in terms of \ $P_n\
.$

\noindent By definition of \ $P_n$ \ we find that:

\vspace{-.3cm}
\begin{itemize}
\item[-] the number of sequences with complexity \ $k+1$ \ and
length \ $n$ \ is \ $\alpha^n P_n (k+1)\ ,$

\vspace{-.3cm}
\item[-] the number of exact sequences with complexity \ $k+1$ \ and
length \ $n$ \ is \ $\alpha^n P_n ((k+1)_e)\ ,$

\vspace{-.3cm}
\item[-] the number of exact sequences with complexity \ $k$ \ and
length \ $n$ \ is \ $\alpha^n P_n (k_e)\ .$
\end{itemize}
\vspace{-.3cm}

Taking into account the definitions of complexity and exact
sequences we conclude that every sequence with complexity \ $k+1$ \
and length \ $n+1$ \ can be obtained from a sequence of length \
$n$ \ in one of the following two ways only:

\vspace{-.3cm}
\begin{itemize}
\item[-] by adding a symbol to a sequence with complexity \ $k+1$ \
which is not exact,

\vspace{-.3cm}
\item[-] by adding a symbol to an exact sequence with complexity \
$k\ .$

\end{itemize}
\vspace{-.3cm}

\noindent We also see that all sequences obtained from exact
sequences of length \ $n$ \ and complexity \ $k+1$ \ by adding a
symbol from \ ${\cal A}$ \ will increase their complexity to \
$k+2$ \ and the number of such sequences is \ $\alpha \cdot
\alpha^n P_n ((k+1)_e)\ .$ From the definition of \ $P_{n+1} (k+1)$
\ and the above observations we conclude that
\bq\label{6}
P_{n+1} (k+1) = {\ds\frac{\alpha \cdot\alpha^n P_n (k+1) - \alpha \cdot
\alpha^n P_n ((k+1)_e) + \alpha \cdot \alpha^n P_n (k_e)}{\alpha^{n+1}}}
\eq
and thus
\bq\label{7}
P_{n+1} (k+1) =  P_n (k+1) + P_n (k_e) - P_n ((k+1)_e)
\eq
for every \ $n, k \in \nn\ .$

Replacing \ $n+1$ \ by \ $n$ \ we have
\bq\label{8}
P_{n} (k+1) =  P_{n-1} (k+1) + P_{n-1} (k_e) - P_{n-1} ((k+1)_e)\ .
\eq
Substituting (8) into (7) we obtain
\bq\label{9}
P_{n+1} (k+1) = P_{n-1} (k+1) + P_{n-1} (k_e) - P_{n-1} ((k+1)_e) +
P_n (k_e)  - P_n ((k+1)_e).
\eq
We replace \ $n$ \ by \ $n-1$ \ in (8) and insert the result in
(9). Continuing this process we arrive at
\bq\label{10}
P_{n+1} (k+1) = P_1 (k+1) + {\ds\sum\limits^n_{r=1}} [P_r (k_e)-P_r
 ((k+1)_e) ]\ .
\eq
Since \ $P_1 (k+1) = 0$ \ for \ $k \geq 1$ \ we have
\bq\label{11}
P_{n+1} (k+1) = {\ds\sum\limits^n_{r=1}} P_r (k_e) -
{\ds\sum\limits^n_{r=1}} P_r  ((k+1)_e )
\eq
for every \ $k,n \in \nn\ .$

\noindent Now, replacing in (11) \ $k$ \ by \ $k+1\ , \ k+2\ , \
k+3,\ldots ,k+(n-k)-1$ \ we obtain
\bq\label{12}
\ba{l}
P_{n+1} (k+2) = {\ds\sum\limits^n_{r=1}} P_r ((k+1)_e ) -
{\ds\sum\limits^n_{r=1}} P_r  ((k+2)_e )\\[2ex]
P_{n+1} (k+3) = {\ds\sum\limits^n_{r=1}} P_r ((k+2)_e ) -
{\ds\sum\limits^n_{r=1}} P_r  ((k+3)_e )\\[2ex]
\hspace{3cm}\vdots\\[2ex]
\hspace{.6cm}P_{n+1} (n) = {\ds\sum\limits^n_{r=1}} P_r ((n-1)_e ) -
{\ds\sum\limits^n_{r=1}} P_r  (n_e)\ .
\ea
\eq
Adding (11) and the above equations and taking into account the
fact that \\ ${\ds\sum\limits^n_{r=1}} P_r (n_e)~=~0$ \ for \ $n \geq
2$ \ we have
\bq\label{13}
{\ds\sum\limits^{n-k}_{s=1}} P_{n+1} (k+s) =
{\ds\sum\limits^n_{r=1}} P_r  (k_e )\ .
\eq
One can easily see that \ $P_{n+1} (k+s)=0$ \ for \ $s>n-k\ ,$
where \ $n>k\geq 1\ .$

Thus, we obtain the following expression for the distribution
function of \ $C_{n+1} :$
\bq\label{14}
P_{n+1}  (C_{n+1} \leq k) = 1 - {\ds\sum\limits^n_{r=1}} P_r  (k_e )\ ,
\eq
which finishes the proof.

 {\it Corollary 1}: For every \ $n$ \ and \ $k\ ,$ \
\bq\label{15}
P_{n+1}  (C_{n+1} \leq k) = P_n (C_n \leq k) -  P_n  (k_e )\ .
\eq

 {\it Proof:} From (14) we have
\bq\label{16}
1 - {\ds\sum\limits^{n-1}_{r=1}} P_r (k_e )= P_n (C_n \leq k)\ .
\eq
Adding (14) and (16) we obtain (15).

 {\it Remark 2}: It follows from the above corollary that
 $ P_{n+1}  (C_{n+1} \leq k) \leq \newline P_n (C_n \leq k).$

 {\it Corollary 2}: From (14) and the fact that [4]
\bq\label{17}
\lim\limits_{n\rightarrow \infty}  P_n (C_n \leq k)= 0
\eq
we deduce that
\bq\label{18}
{\ds\sum\limits^\infty_{r=1}} P_r (k_e )= 1\ .
\eq

\vspace{.5cm}
\centerline{\sc IV. Final Remarks}

The complexity of sequences was suggested
as a statistical test of randomness of a random number generators
and block ciphers [5], [7]. It was
proved in [4] that $ \lim\limits_{n\rightarrow
\infty}P_n(C_n \leq {\ds\frac{n}{\log_\alpha n}})=0.$ Therefore, the sets \
$K_{n,k}: = \{ S \in {\cal A}^n: C_n (S) \leq k\}$ \
seem to be good candidates for critical sets (usually \ $k$ \
is assumed [7] to be \ ${\ds\frac{n}{\log_\alpha n}}).$ This means, in fact,
that for an arbitrarily chosen probability p close to 0 there is $n_0$  such that
for $n>n_0$, for a given randomly chosen sequence S the inequality
$C_n(S)\leq {\ds\frac{n}{\log_\alpha n}}$
holds with probability less than p.
Thus, it is
essential  to estimate \ $ P_n (K_{n,k}) = {\ds\sum\limits^k_{s=1}}
P_n(s)\ ,$ i.e. the levels of significance for \ $ K_{n,k}\ .$ In
practice, for a fixed \ $n$ \ these sums are computed numerically by
finding all terms. Formula (15) makes it possible to find the probability \
$P_{n+1} (K_{n+1,k})$ \ for sequences of length \ $n+1$ \ from the
probabilities \ $P_n(K_{n,k})$ \ and \ $P_n(k_e)$ \ for sequences
of length \ $n$ \ (the latter two can be calculated simultaneously). This
reduces the computation time.\\[2mm]

\newpage
\centerline{\sc References}

\begin{enumerate}

\vspace{-.3cm}
\item[{[1]}] G. Chaitin, "Information-theoretic limitations of formal
systems", J. Ass. Comput. Mach., vol. 21 (1974), pp. 403-424.

\vspace{-.3cm}
\item[{[2]}]
E. Gilbert, T. Kadota, "The Lempel-Ziv algorithm and message
complexity", IEEE Trans. Inform. Theory, vol. 38 (1992), 1839-1842.

\vspace{-.3cm}
\item[{[3]}] A. N. Kolmogorov, "Three approaches to the qualitative
definition of information" Prob. Inform.Transmission, vol. 1
(1965), pp. 1-7.

\vspace{-.3cm}
\item[{[4]}] A. Lempel, J. Ziv, "On the complexity of finite
sequences", IEEE Trans. Inform. Theory, vol. IT-22. No. 1
(1976), pp. 75-81.

\vspace{-.3cm}
\item[{[5]}] A. K. Leung, S. E. Tavares, "Sequence complexity as a
test for cryptographic systems", Advances in Cryptology, Crypto'84,
Springer Verlag LNCS 196  ed. G.R.~Blakley, D.~Chaum (1985),
pp. 468-474.

\vspace{-.3cm}
\item[{[6]}]
G. Louchard, W. Szpankowski, J. Tang, "Average profile of the generalized
digital search tree and the generalized Lempel-Ziv algorithm",
SIAM J. on Compution 28: (3) (1999), pp. 904-934.

\vspace{-.3cm}
\item[{[7]}]
G. Wignarajah, "Complexity tests for statistical independence",
M.S. thesis, University of Toledo, 1985.

\vspace{-.3cm}
\item[{[8]}]
J. Ziv, A. Lempel, "A universal algorithm for sequential date
compression", IEEE Trans. Inform. Theory, vol. 23,
(1977), pp. 337-343.

\vspace{-.3cm}
\item[{[9]}]
J. Ziv, A. Lempel, "Compression of individual sequences via
variable rate coding", IEEE Trans. Inform. Theory, vol. 24,
(1978), pp. 530-536.

\vspace{-.3cm}
\item[{[10]}]
J. Ziv, "Compression, tests for randomness and estimating the
statistical model of individual sequences" in SEQUENCES, R. Capocelli,
Ed. New York: Springer-Verlag 1990, pp. 366-373.

\end{enumerate}

\end{document}